# DETERMINATION OF ONE UNKNOWN THERMAL COEFFICIENT THROUGH THE ONE-PHASE FRACTIONAL LAMÉ-CLAPEYRON-STEFAN PROBLEM

**Domingo Alberto Tarzia**

CONICET - Departamento de Matemática, FCE, Univ. Austral,
Paraguay 1950, S2000FZF  Rosario, Argentina.

*E-mail: [DTarzia@austral.edu.ar](mailto:DTarzia@austral.edu.ar)*

**Abstract:**

We obtain explicit expressions for one unknown thermal coefficient (among the conductivity, mass density, specific heat and latent heat of fusion) of a semi-infinite material through the one-phase fractional Lamé-Clapeyron-Stefan problem with an over-specified boundary condition on the fixed face $x = 0$.

The partial differential equation and one of the conditions on the free boundary include a time Caputo's fractional derivative of order $0 < \alpha < 1$.

Moreover, we obtain the necessary and sufficient conditions on data in order to have a unique solution by using recent results obtained for the fractional diffusion equation exploiting the properties of the Wright and Mainardi functions, given in:

- Roscani - Santillan Marcus, Fract. Calc. Appl. Anal., 16 (2013), 802-815;
- Roscani-Tarzia, Adv. Math. Sci. Appl., 24 (2014), 237-249 and
- Voller, Int. J. Heat Mass Transfer, 74 (2014), 269-277.

This work generalizes the method developed for the determination of unknown thermal coefficients for the classical Lamé-Clapeyron-Stefan problem given in Tarzia, Adv. Appl. Math., 3 (1982), pp. 74-82, which are recovered by taking the limit when the order $\alpha \rightarrow 1^-$.

**2010 AMS Subjects Classification**: 26A33, 35C05, 35R11, 35R35, 80A22.

## 1. INTRODUCTION

Heat transfer problems with a phase-change such as melting and freezing have been studied in the last century due to their wide scientific and technological applications, see [1, 2, 4, 5, 7, 9, 14, 21].

A review of a long bibliography on moving and free boundary problems for phase-change materials (PCM) for the heat equation was given in [24]. The importance of obtaining explicit solutions to some free boundary problems was given in the work [25].

We consider a semi-infinite material, with constant thermal coefficients, which is initially solid at its melting temperature $T_m$. At time $t = 0$, we impose a constant temperature $T_0$ $(T_0 > T_m)$ at the fixed face $x = 0$, and a solidification process begins.

We consider that one of the four thermal coefficients is unknown and that it will be determined by a fractional phase-change problem by imposing an over-specified heat flux condition of the type described in [2, 22, 23].

Fractional differential equations have been developed in the last decades, see for example the books [12,15,18], and the articles [8, 13, 16, 17], and some papers on the fractional Lamé-Clapeyron-Stefan problem were published in the last few years, see [6, 10, 11, 19, 20, 26, 27].

In this paper, the differential equation and a governing condition for the free boundary include a fractional time derivative of order $0 < \alpha < 1$ in the Caputo sense, which is defined as [3]:

$$D^\alpha f(t) = \frac{1}{\Gamma(1-\alpha)} \int_0^t \frac{f^{'}(\tau)}{(t-\tau)^\alpha} \, d\tau \quad \text{if} \quad 0 < \alpha < 1$$

$$= f^{'}(t) \qquad \qquad \text{if} \quad \alpha = 1 \tag{1}$$

where $\Gamma$ is the Gamma Function defined by:

$$\Gamma(x) = \int_0^{+\infty} t^{x-1} \exp(-t) \, dt \, . \tag{2}$$

We also define two very important functions, which will be useful in the next Section:

i) Wright Function [28]:

$$W(z; \alpha, \beta) = \sum_{n=0}^{+\infty} \frac{z^n}{n \, \Gamma(n\alpha + \beta)}, \quad z \in \mathbb{C}, \quad \alpha > -1, \quad \beta \in \mathbb{R} \, . \tag{3}$$

ii) Mainardi Function [8]:

$$M_\upsilon(z) = W(-z; -\upsilon, 1-\upsilon) = \sum_{n=0}^{+\infty} \frac{(-z)^n}{n\,\Gamma(-n\upsilon + 1 - \upsilon)}, \quad z \in \mathbb{C}, \quad \upsilon < 1. \quad (4)$$

We note that the Mainardi function is a particular case of the Wright Function.

Some basic properties for the Caputo fractional derivative and for the Wright function are the following:

$$\frac{\partial W}{\partial z}(z; \alpha, \beta) = W(z; \alpha, \alpha + \beta), \qquad D^\alpha(t^\beta) = \frac{\Gamma(1+\beta)}{\Gamma(1+\beta-\alpha)} t^{\beta-\alpha}, \quad (5)$$

$$W\left(-x; -\frac{1}{2}, 1\right) = \mathrm{erfc}\left(\frac{x}{2}\right), \quad 1 - W\left(-x; -\frac{1}{2}, 1\right) = \mathrm{erf}\left(\frac{x}{2}\right), \quad (6)$$

where the classical error and the complementary error functions are defined by:

$$\mathrm{erf}(x) = \frac{2}{\sqrt{\pi}} \int_0^x e^{-u^2} du, \qquad \mathrm{erfc}(x) = 1 - \mathrm{erf}(x) = \frac{2}{\sqrt{\pi}} \int_x^{+\infty} e^{-u^2} du \quad (7)$$

The method for the determination of unknown thermal coefficients through a one-phase fractional Lamé-Clapeyron-Stefan problem with an over-specified boundary condition at the fixed face $x = 0$ is a new and original problem, which is defined as: Finding the free boundary $x = s(t)$, and the temperature $T = T(x,t)$, and one thermal coefficient such that the following equation and conditions are satisfied:

$$D^\alpha T = \lambda^2 T_{xx}, \qquad 0 < x < s(t), \quad t > 0, \quad (8)$$

$$s(0) = 0, \quad (9)$$

$$T(x,0) = T(+\infty, t) = T_m, \qquad x > 0, \quad t > 0, \quad (10)$$

$$T(s(t), t) = T_m, \qquad t > 0 \qquad , \quad (11)$$

$$-kT_x(s(t), t) = \rho\ell\, D^\alpha s(t), \quad t > 0, \quad (12)$$

$$T(0, t) = T_0 > T_m, \qquad t > 0, \quad (13)$$

$$kT_x(0, t) = -\frac{q_0}{t^{\alpha/2}}, \quad t > 0, \quad (14)$$

where $\rho$ is the density of mass, k is the hermal conductivity, $c$ is the specific heat by unit of mass, $\ell$ is the latent heat of fusion by unit of mass, $\lambda^2 = \dfrac{k}{\rho c} > 0$ is the diffusion coefficient, $T_0\,(> T_m)$ is the temperature at the fixed face $x = 0$ and $q_0 > 0$ is the coefficient which characterized the heat flux at the fixed face $x = 0$. We assume that data $T_0$ and $q_0$ are determined experimentally.

The unknown thermal coefficient can be chosen among the four following ones: $k, \rho, c$ y $\ell$.

The goal of the present work is to obtain in Section II:

a) The solution of the one-phase time fractional Lamé-Clapeyron-Stefan of order $0 < \alpha < 1$ (8)-(13) with an over-specified boundary condition of the heat flux type (14) by giving the explicit expression of the temperature $T = T(x, t)$, the free boundary $x = s(t)$ and the unknown thermal coefficient for the four different cases (see Table 1);

b) The restrictions on the data of the corresponding problem for the four different cases in order to have a unique explicit solution (see Table 1).

We remark that the results and explicit formulae obtained in [23] for the determination of one unknown thermal coefficient through the classical one-phase Lamé-Clapeyron-Stefan problem are generalized for the fractional case $0 < \alpha < 1$, and they can be recovered when $\alpha \to 1^-$ (see Table 2).

## II. DETERMINATION OF ONE UNKNOWN THERMAL COEFFICIENT

First, we obtain a preliminary property in order to have a solution to problem (8)-(14).

**Lemma 1** *The solution of the problem (8)-(14) with $0 < \alpha < 1$ and one unknown thermal coefficient is given by:*

$$s(t) = \lambda \, \xi \, t^{\alpha/2}, \quad \xi > 0 \,, \tag{15}$$

$$T(x, t) = T_0 - \frac{T_0 - T_m}{1 - W\left(-\xi; -\frac{\alpha}{2}, 1\right)}\left[1 - W\left(-\frac{x}{\lambda \, t^{\alpha/2}}; -\frac{\alpha}{2}, 1\right)\right], \tag{16}$$

*where the dimensionaless coefficient $\xi > 0$ and the unknown thermal coefficient must satisfy the following system of equations:*

$$\frac{k(T_0 - T_m)}{\lambda \, q_0 \, \Gamma\left(1 - \frac{\alpha}{2}\right)} = 1 - W\left(-\xi; -\frac{\alpha}{2}, 1\right), \tag{17}$$

$$\frac{c(T_0 - T_m)\Gamma\left(1 - \frac{\alpha}{2}\right)}{\ell \, \Gamma\left(1 + \frac{\alpha}{2}\right)} = \xi \, \frac{1 - W\left(-\xi; -\frac{\alpha}{2}, 1\right)}{M_{\alpha/2}(\xi)}. \tag{18}$$

**Proof.** Following [19, 20], by using properties (5) and (7) we have that the expressions (16) and (15) for the temperature and the free boundary satisfy equation (8) and conditions (9), (10), (11) and (13). Exploiting conditions (12) and (14), we obtain that the dimensionaless coefficient $\xi$ and the unknown thermal coefficient must satisfy conditions (17) and (18). □

Now, we will study the four following cases:

Case 1: Determination of $\{\xi, c\}$ ;

Case 2: Determination of $\{\xi, \ell\}$ ;

Case 3: Determination of $\{\xi, k\}$ ;

Case 4: Determination of $\{\xi, \rho\}$ ,

whose results are summarized in Table 1.

**Remark 1**

*In a analogous manner, we can compute the explicit formulae for the four thermal coefficients of the solid phase of the semi-infinite material by using a solidification process instead of a fusion process.*

**Theorem 2 (Case 1: Determination of the thermal coefficient $c$)**

*If data verify the condition:*

$$\frac{k\rho\ell\left(T_0 - T_m\right)\Gamma\left(1 + \dfrac{\alpha}{2}\right)}{q_0^2\,\Gamma\left(1 - \dfrac{\alpha}{2}\right)} < 1, \tag{19}$$

*then the solution of the case 1 (problem* (8)-(14) *with* $0 < \alpha < 1$ *and the unknown thermal coefficient $c$) is given by:*

$$c = c_\alpha = \frac{\ell\,\Gamma\left(1 + \dfrac{\alpha}{2}\right)}{\left(T_0 - T_m\right)\Gamma\left(1 - \dfrac{\alpha}{2}\right)}\frac{\xi_\alpha}{M_{\alpha/2}(\xi_\alpha)}\left[1 - W\left(-\xi_\alpha; -\frac{\alpha}{2}, 1\right)\right], \tag{20}$$

*where the coefficient* $\xi = \xi_\alpha = \xi_\alpha(\alpha, k, \rho, \ell, q_0, T_0 - T_m) > 0$ *is the unique solution of the equation:*

$$\frac{1 - W\left(-x; -\dfrac{\alpha}{2}, 1\right)}{x} = \frac{k\rho\ell\left(T_0 - T_m\right)\Gamma\left(1 + \dfrac{\alpha}{2}\right)}{q_0^2\,\Gamma^3\left(1 - \dfrac{\alpha}{2}\right)}\frac{1}{M_{\alpha/2}(x)}, \quad x > 0. \tag{21}$$

*Moreover, the temperature* $T(x,t) = T_\alpha(x,t)$ *and the free boundary* $s(t) = s_\alpha(t)$ *are given by the following expressions* $(0 < \alpha < 1)$

$$T_\alpha(x,t) = T_0 - \frac{T_0 - T_m}{1 - W\left(-\xi_\alpha; -\dfrac{\alpha}{2}, 1\right)}\left[1 - W\left(-\frac{x}{\lambda_\alpha\,t^{\alpha/2}}; -\frac{\alpha}{2}, 1\right)\right], \tag{22}$$

$$s_\alpha(t) = \lambda_\alpha\,\xi_\alpha\,t^{\alpha/2}, \tag{23}$$

*the dimensionless coefficient* $\xi_\alpha = \xi_\alpha(\alpha, k, \rho, \ell, q_0, T_0 - T_m) > 0$ *is the unique solution of the equation* (21) *and the diffusion coefficient* $\lambda_\alpha^2$ *is given by the following expression:*

$$\lambda_\alpha^2 = \frac{k}{\rho c_\alpha} = \frac{k(T_0 - T_m)\Gamma\left(1 - \dfrac{\alpha}{2}\right)}{\rho \ell \Gamma\left(1 + \dfrac{\alpha}{2}\right)} \frac{M_{\alpha/2}(\xi_\alpha)}{\xi_\alpha\left[1 - W\left(-\xi_\alpha; -\dfrac{\alpha}{2}, 1\right)\right]} \quad . \tag{24}$$

**Proof.** From condition (18) we obtain expression (20). Taking into account the definition of the diffusion coefficient and the expression (20), from condition (17) we obtain the equation (21) for the dimensionless coefficient $\xi_\alpha$. The equation (21) has a unique positive solution if and only if data verify condition (19). In order to prove this fact we can see that the following real functions

$$G_{1\alpha}(x) = M_{\frac{\alpha}{2}}(x), \quad G_{2\alpha}(x) = W\left(-x; -\frac{\alpha}{2}, 1\right), \quad G_{3\alpha}(x) = 1 - W\left(-x; -\frac{\alpha}{2}, 1\right), \quad x > 0, \tag{25}$$

have the following properties [19]:

$$G_{1\alpha}(0^+) = \frac{1}{\Gamma\left(1 - \dfrac{\alpha}{2}\right)} > 0, \quad G_{1\alpha}(+\infty) = 0, \quad G_{1\alpha}'(x) < 0, \quad \forall x > 0, \tag{26}$$

$$G_{2\alpha}(0^+) = 1, \quad G_{2\alpha}(+\infty) = 0, \quad G_{2\alpha}'(x) < 0, \quad \forall x > 0, \tag{27}$$

$$G_{3\alpha}(0^+) = 0, \quad G_{3\alpha}(+\infty) = 1, \quad G_{3\alpha}'(x) > 0, \quad \forall x > 0, \tag{28}$$

and the real function

$$F_{5\alpha}(x) = \frac{1 - W\left(-x; -\dfrac{\alpha}{2}, 1\right)}{x}, \quad x > 0, \tag{29}$$

is a positive strictly decreasing function because

$$F_{5\alpha}(0^+) = \frac{1}{\Gamma\left(1 - \dfrac{\alpha}{2}\right)} > 0, \quad F_{5\alpha}(+\infty) = 0, \tag{30}$$

and

$$F_{5\alpha}'(x) = \frac{x G_{3\alpha}'(x) - G_{3\alpha}(x)}{x^2} = \frac{x M_{\frac{\alpha}{2}}(x) - G_{3\alpha}(x)}{x^2} < 0, \quad \forall x > 0, \tag{31}$$

owing to the fact

$$G_{3\alpha}(x) = G_{3\alpha}(x) - G_{3\alpha}(0^+) = \int_0^x G_{3\alpha}'(t)\,dt = \int_0^x M_{\frac{\alpha}{2}}(t)\,dt > \int_0^x M_{\frac{\alpha}{2}}(x)\,dt = x M_{\frac{\alpha}{2}}(x), \quad \forall x > 0. \tag{32}$$

Then, we get the expression (24) for the diffusion coefficient. □

**Theorem 3** *If the parameter $\alpha \to 1^-$ then, under the hypothesis* (19), *the solution of the case 1, given by* (22), (23), (20) *and* (21) *coincides with the one given in* [23]:

$$s_1(t) = 2\lambda_1\, \mu_1\, t^{1/2}, \quad \mu_1 > 0\,, \tag{33}$$

$$T_1(x,t) = T_0 - \frac{T_0 - T_m}{erf(\mu_1)}\, erf\left(\frac{x}{2\lambda_1 t^{1/2}}\right), \tag{34}$$

$$c_1 = \frac{\pi\, q_0^2}{k\rho\,(T_0 - T_m)^2}\, erf^2(\mu_1), \qquad \lambda_1 = \sqrt{\frac{k}{\rho c_1}} = \frac{k\,(T_0 - T_m)}{q_0 \sqrt{\pi}\; erf(\mu_1)}\,, \tag{35}$$

*where the dimensionless coefficient $\mu_1 > 0$ is the unique solution of the equation:*

$$\exp(x^2) = \frac{q_0^2 \sqrt{\pi}}{k\,\rho\,\ell\,(T_0 - T_m)}\, \frac{erf(x)}{x}, \quad x > 0\,. \tag{36}$$

*In particular, the inequality* (19) *is transformed in the following one:*

$$\frac{k\,\rho\,\ell\,(T_0 - T_m)}{2\,q_0^2} < 1\,. \tag{37}$$

**Proof.** It follows from (6) and properties of functions $\Gamma$ and $W$. □

**Theorem 4 (Case 2: Determination of the thermal coefficient $\ell$ )**

*If data verify the condition:*

$$\frac{\sqrt{\rho c k}\,(T_0 - T_m)}{q_0\, \Gamma\left(1 - \dfrac{\alpha}{2}\right)} < 1\,, \tag{38}$$

*then the solution of the case 2 (problem* (8)-(14) *with $0 < \alpha < 1$ and the unknown thermal coefficient $\ell$ ) is given by:*

$$\ell = \ell_\alpha = \frac{c\,(T_0 - T_m)\,\Gamma\left(1 - \dfrac{\alpha}{2}\right)}{\Gamma\left(1 + \dfrac{\alpha}{2}\right) F_{4\alpha}(\xi)}\,, \tag{39}$$

*where the coefficient $\xi_\alpha = \xi_\alpha(\alpha, k, \rho, c, q_0, T_0 - T_m) > 0$ is the unique solution of the equation:*

$$1 - W\left(-x; -\frac{\alpha}{2}, 1\right) = \frac{\sqrt{\rho c k}\,(T_0 - T_m)}{q_0\,\Gamma\left(1 - \frac{\alpha}{2}\right)}, \quad x > 0. \tag{40}$$

*Moreover, the temperature* $T(x,t) = T_\alpha(x,t)$ *and the free boundary* $s(t) = s_\alpha(t)$ *are given by the following expressions* $(0 < \alpha < 1)$:

$$T_\alpha(x,t) = T_0 - \frac{T_0 - T_m}{1 - W\left(-\xi_\alpha; -\frac{\alpha}{2}, 1\right)} \left[1 - W\left(-\frac{x}{\lambda_\alpha\, t^{\alpha/2}}; -\frac{\alpha}{2}, 1\right)\right], \tag{41}$$

$$s_\alpha(t) = \lambda\,\xi_\alpha\, t^{\alpha/2}, \tag{42}$$

*the dimensionless coefficient* $\xi = \xi_\alpha = \xi_\alpha(\alpha, k, \rho, c, q_0, T_0 - T_m) > 0$ *is the unique solution of the equation* (40) *and the diffusion coefficient* $\lambda_\alpha^2$ *is given by the following expression:*

$$\lambda_\alpha^2 = \lambda^2 = \frac{k}{\rho c} \quad . \tag{43}$$

**Proof.** From (17) we obtain the equation (40) for the coefficient $\xi_\alpha$, which has a unique solution if and only if data verify the condition (38) because function $G_{3\alpha}$ is a positive increasing function that satisfies the properties given in (28). From (18) we obtain the expression (39), and then we have that (43) holds. □

**Theorem 5** *If the parameter* $\alpha \to 1^-$ *then, under the hypothesis* (38), *the solution of the case 2 given by* (41), (42), (39) *and* (40) *coincides with the one given in* [23]:

$$s_1(t) = 2\lambda\,\mu_1\, t^{1/2}, \quad \mu_1 > 0 \; , \tag{44}$$

$$T_1(x,t) = T_0 - q_0 \sqrt{\frac{\pi}{\rho c k}}\, erf\left(\frac{x}{2\lambda t^{1/2}}\right), \tag{45}$$

$$\ell_1 = q_0 \sqrt{\frac{c}{\rho k}}\, \frac{\exp(-\mu_1^2)}{\mu_1}, \qquad \lambda = \sqrt{\frac{k}{\rho c}} \; , \tag{46}$$

*where the dimensionless coefficient* $\mu_1 > 0$ *is the unique solution of the equation:*

$$erf(x) = \frac{(T_0 - T_m)}{q_0} \sqrt{\frac{\rho c k}{\pi}}, \quad x > 0. \tag{47}$$

*In particular, the inequality* (38) *is transformed in the following one:*

$$\frac{(T_0 - T_m)}{q_0} \sqrt{\frac{\rho c k}{\pi}} < 1. \tag{48}$$

**Proof.** It follows from (6) and properties of functions $\Gamma$ and $W$.  □

### Theorem 6 (Case 3: Determination of the thermal coefficient *k*)

*For any data, the solution of the case 3 (problem (8)-(14) with $0 < \alpha < 1$ and the unknown thermal coefficient $k$) is given by:*

$$k = k_\alpha = \frac{q_0^2\, \Gamma^2\left(1 - \frac{\alpha}{2}\right)}{\rho c\left(T_0 - T_m\right)^2}\left[1 - W\left(-\xi_\alpha; -\frac{\alpha}{2}, 1\right)\right]^2, \qquad (49)$$

*where the coefficient $\xi_\alpha = \xi_\alpha(\alpha, \ell, \rho, c, q_0, T_0 - T_m) > 0$ is the unique solution of the equation:*

$$\frac{x\left[1 - W\left(-x; -\frac{\alpha}{2}, 1\right)\right]}{M_{\alpha/2}(x)} = \frac{c\left(T_0 - T_m\right)\Gamma\left(1 - \frac{\alpha}{2}\right)}{\ell\, \Gamma\left(1 + \frac{\alpha}{2}\right)}, \quad x > 0. \qquad (50)$$

*Moreover, the temperature $T(x,t) = T_\alpha(x,t)$ and the free boundary $s(t) = s_\alpha(t)$ are given by the following expressions $(0 < \alpha < 1)$:*

$$T_\alpha(x,t) = T_0 - \frac{T_0 - T_m}{1 - W\left(-\xi_\alpha; -\frac{\alpha}{2}, 1\right)}\left[1 - W\left(-\frac{x}{\lambda_\alpha\, t^{\alpha/2}}; -\frac{\alpha}{2}, 1\right)\right], \qquad (51)$$

$$s_\alpha(t) = \lambda_\alpha\, \xi_\alpha\, t^{\alpha/2}, \qquad (52)$$

*the dimensionless coefficient $\xi_\alpha = \xi_\alpha(\alpha, c, \rho, \ell, q_0, T_0 - T_m) > 0$ is the unique solution of the equation (50) and the diffusion coefficient $\lambda_\alpha^2$ is given by the following expression:*

$$\lambda_\alpha = \frac{q_0\, \Gamma\left(1 - \frac{\alpha}{2}\right)}{\rho c\left(T_0 - T_m\right)}\left[1 - W\left(-\xi_\alpha; -\frac{\alpha}{2}, 1\right)\right]\ . \qquad (53)$$

**Proof.** From (18) we have that the coefficient $\xi_\alpha$ satisfies the equation (50), which has a unique solution for any data since the real function

$$F_{4\alpha}(x) = \frac{x\left[1 - W\left(-x; -\frac{\alpha}{2}, 1\right)\right]}{M_{\alpha/2}(x)}, \quad x > 0\ . \qquad (54)$$

is a positive increasing function since the following properties hold [19,20]:

$$F_{4\alpha}(0^+) = 0, \quad F_{4\alpha}(+\infty) = +\infty, \quad F'_{4\alpha}(x) > 0, \quad \forall x > 0, \tag{55}$$

Therefore, from (17) we obtain the expressions (49) and (54) for the conductivity $k_\alpha$ and the diffusion coefficient $\lambda_\alpha$, respectively. □

**Theorem 7** *For any data, if the parameter* $\alpha \to 1^-$ *then the solution of the case 3 given by* (51), (52), (49) *and* (50) *coincides with the one given in* [23]:

$$s_1(t) = 2\lambda_1 \mu_1 t^{1/2}, \quad \mu_1 > 0, \tag{56}$$

$$T_1(x,t) = T_0 - \frac{T_0 - T_m}{erf(\mu_1)} erf\left(\frac{x}{2\lambda_1 t^{1/2}}\right), \tag{57}$$

$$k_1 = \frac{\pi q_0^2}{\rho c (T_0 - T_m)^2} \, erf^2(\mu_1), \qquad \lambda_1 = \frac{q_0}{\rho \ell} \frac{\exp(-\mu_1^2)}{\mu_1}, \tag{58}$$

*where the dimensionless coefficient* $\mu_1 > 0$ *is the unique solution of the equation:*

$$x \exp(x^2) \, erf(x) = \frac{c (T_0 - T_m)}{\ell \sqrt{\pi}}, \quad x > 0. \tag{59}$$

**Proof.** It follows from (6) and properties of functions $\Gamma$ and $W$. □

**Theorem 8 (Case 4: Determination of the thermal coefficient** $\rho$ **)**

*For any data, the solution of the case 4 (problem* (8)-(14) *with* $0 < \alpha < 1$ *and the unknown thermal coefficient* $\rho$ *) is given by:*

$$\rho = \rho_\alpha = \frac{q_0^2 \, \Gamma^2\left(1 - \dfrac{\alpha}{2}\right)}{kc(T_0 - T_m)^2} \left[1 - W\left(-\xi_\alpha; -\dfrac{\alpha}{2}, 1\right)\right]^2, \tag{60}$$

*where the coefficient* $\xi_\alpha = \xi_\alpha(\alpha, \ell, k, c, q_0, T_0 - T_m) > 0$ *is the unique solution of the equation* (50). *Moreover, the temperature* $T(x,t) = T_\alpha(x,t)$ *and the free boundary* $s(t) = s_\alpha(t)$ *are given by the following expressions* $(0 < \alpha < 1)$:

$$T_\alpha(x,t) = T_0 - \frac{T_0 - T_m}{1 - W\left(-\xi_\alpha; -\dfrac{\alpha}{2}, 1\right)}\left[1 - W\left(-\frac{x}{\lambda_\alpha\, t^{\alpha/2}}; -\frac{\alpha}{2}, 1\right)\right], \qquad (61)$$

$$s_\alpha(t) = \lambda_\alpha\, \xi_\alpha\, t^{\alpha/2}, \qquad (62)$$

*the dimensionless coefficient $\xi = \xi_\alpha = \xi_\alpha(\alpha, c, k, \ell, q_0, T_0 - T_m) > 0$ is the unique solution of the equation* (50) *and the diffusion coefficient $\lambda_\alpha^2$ is given by the following expression:*

$$\lambda_\alpha = \frac{k\left(T_0 - T_m\right)}{q_0\, \Gamma\left(1 - \dfrac{\alpha}{2}\right)\left[1 - W\left(-\xi_\alpha; -\dfrac{\alpha}{2}, 1\right)\right]} \qquad (63)$$

**Proof.** It is similar to the proof of the case 3 (see Theorem 6). □

**Theorem 9** *For any data, if the parameter $\alpha \to 1^-$ then the solution of the case 4 given by* (61), (62), (60) *and* (50) *coincides with the one given in* [23]:

$$s_1(t) = 2\lambda_1\, \mu_1\, t^{1/2}, \quad \mu_1 > 0\,, \qquad (64)$$

$$T_1(x,t) = T_0 - \frac{T_0 - T_m}{erf(\mu_1)}\, erf\left(\frac{x}{2\lambda_1 t^{1/2}}\right), \qquad (65)$$

$$\rho_1 = \frac{\pi q_0^2}{kc\left(T_0 - T_m\right)^2}\, erf^2(\mu_1), \qquad \lambda_1 = \frac{k\left(T_0 - T_m\right)}{q_0\sqrt{\pi}\, erf(\mu_1)}, \qquad (66)$$

*where the dimensionless coefficient $\mu_1 > 0$ is the unique solution of the equation:*

$$x \exp(x^2)\, erf(x) = \frac{c\left(T_0 - T_m\right)}{\ell\,\sqrt{\pi}}, \quad x > 0\,. \qquad (67)$$

**Proof.** It is similar to the proof of the case 3 (see Theorem 7). □

Now, in order to summarize our results on the determination of one unknown thermal coefficient through a fractional Lamé-Clapeyron-Stefan problem with an over-specified heat flux boundary condition on the fixed face, we show the formulae and restrictions for data for the four cases for the fractional Lamé-Clapeyron-Stefan problem with $0 < \alpha < 1$ (see Table 1) and for the classical Lamé-Clapeyron-Stefan problem with $\alpha = 1$ (see Table 2).

| Case # | Explicit formulae for the unknown thermal coefficient | Equation that must satisfy the parameter $\xi$ | Restrictions on data |
|---|---|---|---|
| 1 | $c = \dfrac{\ell\,\Gamma\left(1+\frac{\alpha}{2}\right)}{(T_0-T_m)\Gamma\left(1-\frac{\alpha}{2}\right)}F_{4\alpha}(\xi)$ | $\dfrac{1-W\left(-x;-\frac{\alpha}{2},1\right)}{x} = \dfrac{k\rho\ell(T_0-T_m)\Gamma\left(1+\frac{\alpha}{2}\right)}{q_0^2\,\Gamma^3\left(1-\frac{\alpha}{2}\right)}\dfrac{1}{M_{\alpha/2}(x)}, \quad x>0$ | $\dfrac{k\rho\ell(T_0-T_m)\Gamma\left(1+\frac{\alpha}{2}\right)}{q_0^2\,\Gamma\left(1-\frac{\alpha}{2}\right)}<1$ |
| 2 | $\ell = \dfrac{c(T_0-T_m)\Gamma\left(1-\frac{\alpha}{2}\right)}{\Gamma\left(1+\frac{\alpha}{2}\right)F_{4\alpha}(\xi)}$ | $G_{5\alpha}(x) = \dfrac{\sqrt{\rho c k}\,(T_0-T_m)}{q_0\,\Gamma\left(1-\frac{\alpha}{2}\right)}, \quad x>0$ | $\dfrac{\sqrt{\rho c k}\,(T_0-T_m)}{q_0\,\Gamma\left(1-\frac{\alpha}{2}\right)}<1$ |
| 3 | $k = \dfrac{q_0^2\,\Gamma^2\left(1-\frac{\alpha}{2}\right)}{\rho c(T_0-T_m)^2}\left[1-W\left(-\xi;-\frac{\alpha}{2},1\right)\right]^2$ | $F_{4\alpha}(x) = \dfrac{c(T_0-T_m)\Gamma\left(1-\frac{\alpha}{2}\right)}{\ell\,\Gamma\left(1+\frac{\alpha}{2}\right)}, \quad x>0$ | ----------- |
| 4 | $\rho = \dfrac{q_0^2\,\Gamma^2\left(1-\frac{\alpha}{2}\right)}{kc(T_0-T_m)^2}\left[1-W\left(-\xi;-\frac{\alpha}{2},1\right)\right]^2$ | $F_{4\alpha}(x) = \dfrac{c(T_0-T_m)\Gamma\left(1-\frac{\alpha}{2}\right)}{\ell\,\Gamma\left(1+\frac{\alpha}{2}\right)}, \quad x>0$ | ----------- |

**Table 1**. Summary of the results corresponding to the determination of one unknown thermal coefficient through a fractional Lamé-Clapeyron-Stefan problem with an over-specified heat flux boundary condition on the fixed face (4 cases).

| Case # | Explicit formulae for the unknown thermal coefficient | Equation that must satisfy the parameter $\mu$ | Restrictions on data |
|---|---|---|---|
| 1 | $c = \dfrac{\pi q_0^2}{\rho k(T_0-T_m)^2}\,erf^2(\mu)$ | $\dfrac{erf(x)}{x} = \dfrac{k\rho\ell(T_0-T_m)}{q_0^2\sqrt{\pi}}\exp(x^2), \quad x>0$ | $\dfrac{\rho\ell k(T_0-T_m)}{2q_0^2}<1$ |
| 2 | $\ell = q_0\sqrt{\dfrac{c}{\rho k}}\,\dfrac{\exp(-\mu^2)}{\mu}$ | $erf(x) = \dfrac{(T_0-T_m)}{q_0}\sqrt{\dfrac{\rho c k}{\pi}}, \quad x>0$ | $\dfrac{(T_0-T_m)}{q_0}\sqrt{\dfrac{\rho c k}{\pi}}<1$ |
| 3 | $k = \dfrac{\pi q_0^2}{\rho c(T_0-T_m)^2}\,erf^2(\mu)$ | $x\exp(x^2)\,erf(x) = \dfrac{c(T_0-T_m)}{\ell\sqrt{\pi}}, \quad x>0$ | ----------- |
| 4 | $\rho = \dfrac{\pi q_0^2}{kc(T_0-T_m)^2}\,erf^2(\mu)$ | $x\exp(x^2)\,erf(x) = \dfrac{c(T_0-T_m)}{\ell\sqrt{\pi}}, \quad x>0$ | ----------- |

**Tabla 2**. Summary of the results corresponding to the determination of one unknown thermal coefficient through a classical Lamé-Clapeyron-Stefan problem ($\alpha=1$) with an over-specified heat flux boundary condition on the fixed face (4 cases). These results were obtained by taking $\alpha\to1^-$ in the results given in Table 1 [23].

### ACKNOWLEDGEMENTS

The present work has been sponsored by the Projects PIP N° 0534 from CONICET – Univ. Austral, and by AFOSR-SOARD Grant FA9550-14-1-0122.